\documentclass[11pt]{article}
\title%
{\vskip-2.0em Splitting maps and norm bounds for the cyclic cohomology of biflat Banach algebras}

\author{Yemon Choi}

\newcommand{\update}{20th April 2010}
\date{\update%
}

\usepackage{amsmath}
\usepackage{amssymb,amsfonts}
\usepackage[usenames]{color}
\usepackage{hyperref}

\newcounter{pulse}[section]
\numberwithin{pulse}{section}  

\RequirePackage[PS,nohug]{diagrams}
\newarrow{Dots} ....>
\newcommand{\Rl}[2]{\pile{\rTo^{#1}\\\lDots_{#2}}}
\newcommand{\ackndiag}{Uses Paul Taylor's {\tt diagrams.sty} macros.}

\usepackage{amsthm}
\newtheorem{thm}[pulse]{Theorem}
\newtheorem{prop}[pulse]{Proposition}
\newtheorem{lem}[pulse]{Lemma}

\theoremstyle{definition}
\newtheorem{dfn}[pulse]{Definition}

\theoremstyle{remark}
\newtheorem{rem}[pulse]{Remark}

\numberwithin{equation}{section}
\numberwithin{figure}{section}

\newcommand{\dt}[1]{\textcolor{Maroon}{\textsf{#1}}}

\newenvironment{newnum}{%
\begin{enumerate}

}{\end{enumerate}\ignorespacesafterend}

\newcommand{\norm}[1]{\Vert{#1}\Vert}

\newcommand{\tp}{\mathop{\otimes}}
\newcommand{\ptp}{\mathop{\widehat{\otimes}}}
\newcommand{\blob}{\bullet}
\newcommand{\sid}{\operatorname{id}}

\newcommand{\defeq}{\overset{\rm def}{=}}
\newcommand{\blank}{\mathop{\underline{\quad}}}

\newenvironment{YCnum}{%
\begin{enumerate}

}{\end{enumerate}\ignorespacesafterend}

\newcommand{\dif}{\delta}

\newcommand{\cCo}[2]{{\mathcal#1}_\lambda^{#2}}

\newcommand{\image}{\mathop{\rm im}}
\newcommand{\inc}{\mathop{\rm inc}}

\newcommand{\cA}{{\mathcal A}}
\newcommand{\cC}{{\mathcal C}}
\newcommand{\cH}{{\mathcal H}}
\newcommand{\cZ}{{\mathcal Z}}

\newcommand{\bt}{\mathop{\textbf{\textit{t}\/}}\nolimits}

\newcommand{\sD}{{\mathsf D}}
\newcommand{\sE}{{\mathsf E}}
\newcommand{\sF}{{\mathsf F}}
\newcommand{\sG}{{\mathsf G}}

\newcommand{\al}{\alpha}
\newcommand{\lm}{\lambda}
\newcommand{\sig}{\sigma}

\begin{document}
\maketitle

\begin{abstract}
We revisit the old result that biflat Banach algebras have the same cyclic cohomology as ${\mathbb C}$, and obtain a quantitative variant (which is needed in separate, joint work of the author on the simplicial and cyclic cohomology of band semigroup algebras).
Our approach does not rely on the Connes-Tsygan exact sequence, but is motivated strongly by its construction as found in \cite{Connes_IHES_NCDG} and \cite{Hel_CT}.
\end{abstract}

\begin{section}{Introduction}
Cyclic cohomology and (simplicial) Hochschild cohomology of an ass\-oc\-iative algebra are related by a long exact sequence, discovered independently by A.~Connes and B.~L. Tsygan and bearing their names. The construction of the Connes-Tsygan exact sequence in the setting of Banach algebras, not necessarily with identity, can be found in~\cite{Hel_CT}; that article contains, among other results, a determination up to topological isomorphism of the continuous cyclic cohomology groups of biflat Banach algebras. (From here onwards, we shall for sake of brevity omit explicit mention of the adjective `continuous'; it is always understood that all cochains being considered on a Banach algebra are assumed to be continuous.)

Recent work of the author with F. Gourdeau and M. C. White~\cite{CGW_bands} calculates the cyclic and simplicial cohomology of certain kinds of Banach algebra~$\cA$: these admit an $\ell^1$-direct sum decomposition at the level of Banach spaces
$$ \cA = \ell^1{\rm-}\bigoplus_{\alpha \in L} A_\al $$
where the indexing set $L$ is a semilattice, each $A_\al$ is a closed, \emph{biflat} subalgebra, and where $A_\al\cdot A_\gamma \subseteq A_{\al\gamma}$ for all $\alpha,\gamma\in L$. Moreover, each $A_\al$ is biflat with constant~$1$ (see Definition~\ref{dfn:biflatness-constant} below).
The arguments in \cite{CGW_bands} are not short: they involve careful, direct calculations with the cyclic cochain complex, and ``induction by hand'' rather than invoking the machinery of spectral sequences.

A~necessary preliminary step in \cite{CGW_bands} is to reduce the problem to the case of those cyclic cocycles $\psi\in \cCo{Z}{n}(\cA)$ which furthermore satisfy the following normalization condition: $\psi(a_1,\dots, a_n)(a_0)=0$ whenever $a_0,\dots, a_n$ all belong to a common component algebra $A_\gamma$, for some $\gamma\in L$.
Since the indexing semilattice $L$ is infinite, in order to carry out this reduction, we need to know not only that `most' cyclic cocycles on each $A_\alpha$ are cyclic coboundaries, but that we can solve the cohomology problem with control of the constants of cobounding. It is insufficient to {merely cite or quote} the known result that $\cCo{H}{n}(A_\alpha)$ is zero in odd dimensions and isomorphic to $\cCo{Z}{0}(A_\alpha)$ in even dimensions; we have to say something about the `constant of openness' of the Hochschild coboundary operator as a bounded linear map from $\cCo{C}{n-1}(A_\alpha)$ to $\cCo{C}{n}(A_\alpha)$.
\medskip

With these observations in mind, the present article has the following aims:
\begin{itemize}
\item[--] to compute the cohomology of the cyclic cochain complex of a biflat Banach algebra, without recourse to other auxiliary or underlying complexes;
\item[--] to do so via `cobounding maps' which can be written down explicitly or recursively, and whose norms are given explicit bounds in terms of the biflatness constant of the given algebra.
\end{itemize}
It is hoped that the exposition is accessible to those familiar with the rudiments of the homological theory of Banach algebras, as found in~\cite{Hel_HBTA}.

\medskip
To be more precise, we need to recall some terminology.
Given a Banach algebra $A$, let $\pi:A\ptp A\to A$ be the bounded linear map defined by $\pi(a\tp b) = ab$, where $\ptp$ denotes the projective tensor product of Banach spaces.

\begin{dfn}\label{dfn:biflatness-constant}
Let $K\geq 1$. We say that a Banach algebra $A$ is \dt{biflat with constant $K$} if there exists a bounded linear $A$-bimodule map $\rho:(A\ptp A)^*\to A^*$, with $\norm{\rho}\leq K$, such that $\rho\pi^*(\psi)=\psi$ for all $\psi\in A^*$.
\end{dfn}

By \cite[Exercise VII.2.8]{Hel_HBTA}, a Banach algebra is biflat in the homological sense if and only if it is biflat (for some constant $K$) in the sense of this definition.

\begin{thm}\label{t:newNEEDED}
Let $A$ be a Banach algebra which is biflat with constant $K$, and let $m\geq 0$\/.
\begin{YCnum}
\item For every $\psi\in \cZ^{2m+1}_\lm(A)$, there exists $\chi\in \cC^{2m}_\lm(A)$ such that $\psi=\dif\chi$ and
\[ \norm{\chi}\leq 2(m+1)^3K^{4m}\norm{\psi} \,.\]   
\item For every $\psi\in \cZ^{2m+2}_\lm(A)$, there exist $\chi\in \cC^{2m+1}_\lm(A)$
 and $\tau\in \cZ^0_\lm(A)$
 such that
$\psi=\tau^{(2m+2)}+\dif\chi$\/, and which satisfy
\[ \norm{\tau}\leq K^{2m+2}\norm{\psi} \quad,\quad \norm{\chi}\leq 2(m+1)^3K^{4m+2}\norm{\psi}\,.\]
Here $\tau^{(2m)}\in\cCo{Z}{2m}$ is the cyclic cocycle defined by
\[ \tau^{(2m)}(a_1,\dots, a_{2m})(a_0) \defeq \tau(a_0\dotsb a_{2m}).\]
\end{YCnum}
\end{thm}

Theorem~\ref{t:newNEEDED} can be thought of as a
quantitative version of \cite[Theorem 25]{Hel_CT}; our goal in this note is to give a \emph{self-contained proof} of this theorem, spelling out the details explicitly. 
Our approach is direct and does not rely on the results of \cite{Hel_CT}, although it is broadly inspired by the same point of view. Moreover, since we have certain applications in mind, we will impose slightly stronger conditions, and hence obtain slightly stronger results.

Although the proof of Theorem~\ref{t:newNEEDED} does not require us to (re)construct the full Connes-Tsygan exact sequence, for sake of completeness we show in an appendix how once can obtain an analogous long exact sequence by pushing the calculations of \S2 a little further.
The reader may also wish to consult~\cite[Part II]{Connes_IHES_NCDG} for some of the underlying motivation behind various {\it ad hoc}\/ calculations in \S2 and the appendix.
For a more conceptual approach to the cyclic cochain complex and various operations on it, see~\cite{Quillen_IHES88}.

\begin{rem}
Even without the constants, Theorem~\ref{t:newNEEDED} is still slightly more precise than the stated result in \cite{Hel_CT}, because in even degrees it identifies explicit representatives of cyclic cohomology classes, namely those arising from traces.
\end{rem}

\end{section}


\begin{section}{An abstract version}\label{s:abstract}
The key calculations required to prove Theorem~\ref{t:newNEEDED} are a matter of judicious diagram chasing on chain complexes. Therefore we shall work in greater generality than will eventually be needed, so as to emphasise the formal nature of the core calculations, and to clarify the structure of our argument.

\subsection*{Notation and other preliminaries.}
We assume familiarity with the basic definitions of (co)chain complexes, chain maps between them, and (chain) homotopies between chain maps. If
\[ \tag{$\sE^\blob$}
 \begin{diagram}
\dots & \rTo & \sE^{n-1} & \rTo^\dif & \sE^n & \rTo^\dif & \sE^{n+1} & \rTo & \dots
\end{diagram}
\]
is a cochain complex, we write $H^n(\sE^\blob)$ for the $n$th cohomology group of this complex, i.e.
\[ H^n(\sE^\blob) = \frac{\ker( \dif: \sE^n \to \sE^{n+1})}{\image(\dif: \sE^{n-1}\to \sE^n)} \]
To reduce notational clutter, we also adopt the convention of omitting the index on individual components of chain maps: thus, if $M: \sE^\blob \to \sF^\blob$ is a chain map, we will denote the component in degree $n$ by $M:\sE^n \to\sF^n$. 

Throughout, $\sid$ denotes the identity map on a given vector space, or chain complex, or algebra; it will hopefully be clear from context what the domain of $\sid$ is.

\subsection*{The main setup}
Consider an exact sequence
\[ \begin{diagram}
0 & \rTo & \sD & \rTo^{\imath} & \sE & \rTo^M & \sF & \rTo^N & \sG & \rTo & 0
\end{diagram} \]
where each of $\sD$\/, $\sE$\/, $\sF$ and $\sG$ are cochain complexes of vector spaces and linear maps, and $\imath$\/, $M$ and $N$ are all chain maps between these complexes. (Much of what follows, being abstract diagram-chasing, would work just as well in any additive category that has kernels and cokernels.)

Let $\dif$, $\dif'$ and $\dif''$ denote the differentials for the complexes $\sE^\blob$, $\sF^\blob$ and $\sG^\blob$ respectively. To simplify notation slightly, we will regard $\sD^\blob$ as a subcomplex of $\sE^\blob$ and $\imath$ as the inclusion map, so that the differential on $\sD^\blob$ will also be denoted by $\dif$.

\begin{rem}
In our intended application, $\sD^\blob$ and $\sG^\blob$ will both be equal to $\cC^\blob_\lm(A)$, the complex of cyclic cochains on some given algebra~$A$\/. 
\end{rem}

We now assume that, for each $n$, there exist linear maps
\[ \begin{diagram}
\sG^n & \rTo^{j_n} & \sF^n & \rTo^{h_n} & \sE^n & \rTo^{P_n} & \sD^n
\end{diagram} \]
such that
\begin{equation}
Nj_n =\sid \quad;\quad Mh_n+j_n N = \sid \quad;\quad h_nM+ \imath P_n =\sid \quad;\quad P_n\imath =\sid\,.
\end{equation}
We shall omit the indexing suffices on the maps $j$, $h$ and $P$.

The following lemma collects some useful identities for later reference: they are easily proved by diagram-chasing and we omit the details.
\begin{lem}[Useful identities]\label{l:useful1}\
\begin{newnum}
\item $\dif h\dif'M = \dif hM\dif = - \dif \imath P \dif = -\imath\dif P\dif$\/.
\item $\dif'j\dif''N =\dif'j N \dif' = - \dif'Mh\dif' = - M\dif h\dif'$\/.
\end{newnum}
\end{lem}

\begin{subsection}{Constructing $\widetilde{S}$}
The definition is straightforward: for each $n$, define $S^\natural:\sG^n\to\sE^{n+2}$ to be the composite map
\begin{equation}
\begin{diagram}[tight,height=2em]
\sE^{n+2} & \\
\uTo^\dif & & \\
\sE^{n+1} & \lTo^h & \sF^{n+1} \\
 & & \uTo^{\dif'} & \\
 & & \sF^n & \lTo^j & \sG^n
\end{diagram}
\end{equation}
and put $\widetilde{S}=PS^\natural$\/.

\begin{lem}[$\widetilde{S}$ is a chain map]
For each $n$\/, we have $S^\natural \dif'' = \imath \dif\widetilde{S}$\/, and hence $\widetilde{S}\dif''=\dif\widetilde{S}$.
\end{lem}

\begin{rem}
Note that in this lemma, we don't need $\sE^{\blob}$ and $\sF^{\blob}$ to be acyclic.
\end{rem}

\begin{figure}
\label{fig:S-is-a-chain-map}
\caption{A diagram chase}
\[
\begin{diagram}[tight,height=2.5em]
0\to & \sD^{n+2} & \rTo^{\imath} & \sE^{n+2} &  &  &  &   \\
 & \uTo^{\dif} & & \uTo^{\dif} & &      & &   \\
 & \sD^{n+1} & \Rl{\imath}{P} & \sE^{n+1} & \Rl{M}{h} & \sF^{n+1} &  &   \\
 & & & \uTo^{\dif} &   & \uTo^{\dif'} & &  \\
 & & & \sE^n & \Rl{M}{h} & \sF^n & \Rl{N}{j} & \sG^n  \\
 & & &    &   &  \uTo^{\dif'}   &     & \uTo^{\dif''} \\
 & & &    &   & \sF^{n-1}  & \Rl{N}{j} & \sG^{n-1}  & \to 0 \\
\end{diagram}
\]
\end{figure}
\begin{proof}
Consider Figure~\ref{fig:S-is-a-chain-map}, in which all squares consisting of solid arrows commute.
We have
\begin{align*}
S^\natural \dif'' = \dif h\dif'j\dif''
	& = \dif h\dif'j\dif''Nj & \text{(since $Nj=\sid$)} \\
	& = - \dif h M\dif h\dif'j & \text{(by Lemma~\ref{l:useful1}(ii))} \\
	& = \imath \dif P\dif h\dif'j & \text{(by Lemma~\ref{l:useful1}(i))} \\
  & = \imath \dif \widetilde{S}
\end{align*}
as required. The last part is immediate since $P\imath=\sid$.
\end{proof}

Since chain maps between chain complexes induce maps between their homology, $\widetilde{S}$ descends to give a map $S: H^n(\sG^\blob)\to H^n(\sD^{\blob+2}) \equiv H^{n+2}(\sD)$\/.
In Section~\ref{s:appl}, when $\sD^\blob=\sG^\blob$ is the complex of cyclic cochains, $S$ will then (up to a trivial rescaling) be the \emph{shift map} in the Connes-Tsygan exact sequence.
\end{subsection}

\begin{subsection}{A homotopy inverse to $\widetilde{S}$, under certain hypotheses}
We now assume contractibility of both $\sE^\blob$ and $\sF^\blob$, given by contracting homotopies $\sig:\sE^\blob\to\sE^{\blob-1}$ and $\sig':\sF^\blob\to\sF^{\blob-1}$ respectively. Thus $\dif\sig+\sig\dif=\sid$ and $\dif'\sig'+\sig'\dif'=\sid$ in each degree.

Again, it will be convenient to collect some useful identities in a separate lemma, for later reference. We omit the proofs, which consist of easy diagram-chasing.
\begin{lem}[More useful identities]\label{l:useful2}\
\begin{newnum}
\item $\dif'M\sig\imath = M\dif\sig\imath = - M\sig\dif\imath = -M\sig\imath\dif$.
\item $\dif''N\sig'M = N\dif'\sig' M = - N\sig'\dif' M = - N\sig'M\dif$. 
\end{newnum}
\end{lem}

We now define $\widetilde{R}:\sD^{n+2}\to\sG^n$ to be the composite mapping
\begin{equation}
\begin{diagram}[tight,height=2em]
\sD^{n+2} & \rTo^{\imath} & \sE^{n+2} \\
 & & \dTo^{\sig} & \\
 & & \sE^{n+1} & \rTo^M & \sF^{n+1} \\
 & & &  & \dTo^{\sig'} \\
 & & &  & \sF^n & \rTo^N & \sG^n\,.
\end{diagram}
\end{equation}

\begin{lem}
$\widetilde{R}$ defines a chain map $\sD^{\blob+2}\to \sG^\blob$\/.
\end{lem}

\begin{figure}
\label{fig:R-is-a-chain-map}
\caption{More chasing}
\[
\begin{diagram}[tight,height=2.5em]
 \sD^{n+2} & \rTo^{\imath} & \sE^{n+2} &  &  &  &   \\
 \uTo^{\dif} & & \uTo^{\dif}\dDots_{\sig} & &      & &   \\
  \sD^{n+1}  & \rTo^{\imath}  & \sE^{n+1} & \rTo^M & \sF^{n+1} &  &   \\
 & & \uTo^{\dif}\dDots_{\sig} & & \uTo^{\dif'}\dDots_{\sig'} & &  \\
 & &  \sE^n  & \rTo^M & \sF^n & \rTo^N & \sG^n  \\
 & &     &   &  \uTo^{\dif'}\dDots_{\sig'}   &     & \uTo^{\dif''} \\
 & &     &   & \sF^{n-1}  & \rTo^N & \sG^{n-1}    \\
\end{diagram}
\]
\end{figure}
\begin{proof}
It will be useful to consider Figure~\ref{fig:R-is-a-chain-map}, where once again all squares consisting of solid arrows commute. 
We have
\begin{flalign*}
\quad\widetilde{R}\dif
  =  N\sig' M\sig\imath\dif
 & = -N\sig' \dif'M\sig\imath & \text{(by Lemma \ref{l:useful2}(i))}\quad \\
 &  = \dif'' N\sig'M\sig\imath  & \text{(by Lemma \ref{l:useful2}(ii))}\quad \\
 & = \dif'' \widetilde{R} \,,
\end{flalign*}
as required.
\end{proof}

\begin{prop}[$\widetilde{S}\widetilde{R}$ is homotopic to $\sid$]
Let $T^\natural \defeq\dif h \sig'M\sig\imath -\sig\imath: \sD^{n+2} \to \sE^{n+1}$\/. Then
\begin{equation}\label{eq:Squeersy-and-me}
S^\natural \widetilde{R} = \imath +T^\natural\dif +\imath\dif PT^\natural\,,
\end{equation}
and so $\widetilde{S}\widetilde{R}= \sid + PT^\natural\dif + \dif PT^\natural$.
\end{prop}

\begin{proof}
It suffices to prove the identity \eqref{eq:Squeersy-and-me}. We have
\begin{flalign*}
\quad S^\natural\widetilde{R}
  = \dif h\dif' j N\sig'M\sig\imath 
 & = \dif h\dif' (\sid-Mh)\sig'M\sig\imath \\
 & = \dif h \dif'\sig'M\sig\imath - \dif h \dif' M h \sig'M\sig\imath \\
 & = \dif h \dif'\sig'M\sig\imath + \imath\dif P\dif h \sig'M\sig\imath  & \text{(by Lemma~\ref{l:useful1}(i)) }\quad \\
 & = \dif h \dif'\sig'M\sig\imath + \imath\dif P(T^\natural+\sig\imath) .
\end{flalign*}
Hence
\begin{flalign*}
\quad S^\natural\widetilde{R} - \imath\dif PT^\natural 
 & =  \dif h \dif'\sig'M\sig\imath + \imath\dif P\sig\imath \\
 & =  \dif h (\sid-\sig'\dif')M\sig\imath + \dif \imath P \sig\imath \\
 & =  \dif(hM+\imath P)\sig\imath - \dif h \sig'\dif'M\sig\imath \\
 & =  \dif(hM+\imath P)\sig\imath + \dif h \sig'M\sig\imath\dif & & \text{(by Lemma~\ref{l:useful2}(i))} \quad \\
 & =  \dif\sig\imath + \dif h \sig'M\sig\imath\dif \\
 & =  (\sid-\sig\dif)\imath + (T^\natural+\sig\imath)\dif 
 &  =  \imath + T^\natural\dif\,,\hbox{\hskip 3em}
\end{flalign*}
where the last step follows because $-\sig\dif\imath+\sig\imath\dif = -\sig\dif\imath+\sig\dif\imath=0$\/.
\end{proof}

\begin{rem}
In particular, we see that under the hypotheses imposed on $\sE^\blob$ and $\sF^\blob$, the map $S: H^*(\sG^\blob) \to H^*(\sD^{\blob+2})$ is an isomorphism of cohomology groups in each degree\/. If we were to construct the full `SBI sequence', then the same conclusion could be obtained without constructing the map $\widetilde{R}$, by inserting $0$ at appropriate places in the sequence.
In an appendix, we shall show how the SBI sequence can be constructed in our setting.
\end{rem}
\end{subsection}

\end{section}

\begin{section}{Application to cyclic cohomology}\label{s:appl}
Recall that, as in \cite[Part II]{Connes_IHES_NCDG}, we are defining the cyclic cohomology of a (Banach) algebra to be the cohomology of the complex of (continuous) cyclic cochains.
%

\begin{subsection}{Notation and reminders}
Throughout, $A$ is a fixed Banach algebra, and $A^*$ its dual. $\cZ(A^*)$ denotes the space of continuous traces on~$A$\/. 

The following definitions are standard, and we repeat them merely to fix notation. (We are, by a slight abuse of notation, identifying continuous $n$-multilinear maps from a Banach algebra $A$ to its dual, with continuous $n+1$-linear multilinear functionals on $A$\/; this is the approach taken in \cite[Part II]{Connes_IHES_NCDG}, for instance, and simplifies some of the formulas that follow.)

Let $n\geq 0$\/. $\cC^n(A)$ denotes the space of bounded $(n+1)$-multilinear functionals on $A$, while $\dif:\cC^n(A)\to \cC^{n+1}(A)$ is given by the usual Hochschild coboundary operator:
\begin{equation}\label{eq:hoch}
\begin{aligned}
\dif\psi(a_0,a_1,\dots, a_{n+1})
 = & \sum_{j=0}^n (-1)^j \psi(a_0, \dots, a_ja_{j+1},\dots, a_{n+1}) \\
 &+ (-1)^{n+1}\psi(a_{n+1}a_0, a_1,\dots,a_n)\;.
\end{aligned}
\end{equation}
We also need to consider the \emph{truncated} Hochschild coboundary operator $\dif':\cC^n(A)\to \cC^{n+1}(A)$ , which is defined by
\begin{equation}\label{eq:hoch'}
\begin{aligned}
\dif'\psi(a_0, a_1,\dots, a_{n+1})
 = & \sum_{j=0}^n (-1)^j \psi(a_0, \dots, a_ja_{j+1},\dots, a_{n+1}) \,.
\end{aligned}
\end{equation}
The ``signed cyclic shift'' operator $\bt: \cC^n(A)\to \cC^n(A)$ is given by
\[ \bt\psi(a_0,a_1,\dots, a_n) = (-1)^n\psi(a_n,a_0,\dots, a_{n-1}) \,.\]
Note that in degree $n$, the operator $\bt$ is periodic of order $n+1$\/.

Elements invariant under the action of $\bt$ are called \dt{cyclic cochains}, and the space of all cyclic $n$-cochains is denoted by $\cC^n_\lm(A)$. Although $\bt$ is not a chain map, it can be shown\footnotemark\ that
\footnotetext{This can be found in several texts. See, for instance, \cite[Part II, Lemma 3]{Connes_IHES_NCDG}.}
\begin{equation}\label{eq:intertwining}
N\dif' = \dif N
\end{equation}
where $N:\cC^n(A)\to\cC^n_\lm(A)$ is defined to be the \dt{averaging operator}
$(n+1)^{-1}\sum_{k=0}^n \bt^k$\/.
In particular, $\dif(\cC^n_\lm(A))\subseteq \cC^n_\lm(A)$, so that the cyclic cochains form a subcomplex $\cC^\blob_\lm(A)$ of the full Hochschild cochain complex. The cohomology groups of the cyclic cochain complex will be denoted by $\cH_\lm^n(A)$, for $n\geq 0$\/.

\medskip
To treat the Hochschild and cyclic cochain complexes within the framework of Section~\ref{s:abstract}, we make the following definitions.
Abusing notation slightly, we take:
\begin{itemize}
\item $\sD^\blob$ to be the complex
\[ \begin{diagram}[tight,width=2.8em,height=1.9em]
 & \sD^{-2}  & & \sD^{-1} & & \sD^0 & & \sD^1 & &  \sD^2 & \\
\dots\to & 0 & \rTo& \cZ(A^*) & \rTo^{\inc} & \cC^0_\lm(A) & \rTo^\dif & \cC^1_\lm(A) & \rTo^\dif & \cC^2_\lm(A) & \rTo^\dif &\dots
\end{diagram} \]
\item $\sE^\blob$ to be the complex
\[ \begin{diagram}[tight,width=2.8em,height=1.9em]
 & \sE^{-2}  & & \sE^{-1} & & \sE^0 & & \sE^1 & &  \sE^2 & \\
\dots\to & 0 & \rTo& \cZ(A^*) & \rTo^{\inc} & \cC^0(A) & \rTo^\dif & \cC^1(A) & \rTo^\dif & \cC^2(A) & \rTo^\dif &\dots
T\end{diagram} \]
\item $\sF^\blob$ to be the complex
\[ \begin{diagram}[tight,width=2.8em,height=1.9em]
 & \sF^{-2}  & & \sF^{-1} & & \sF^0 & & \sF^1 & & \sF^2 & \\
\dots \to & 0 & \rTo& 0 & \rTo & \cC^0(A) & \rTo^{\dif'} & \cC^1(A) & \rTo^{\dif'} & \cC^2(A) & \rTo^{\dif'} & \dots
\end{diagram} \]
\item $\sG^\blob$ to be the complex
\[ \begin{diagram}[tight,width=2.8em,height=1.9em]
 & \sG^{-2}  & & \sG^{-1} & & \sG^0 & & \sG^1 & & \sG^2 & \\
 \dots \to & 0 & \rTo& 0 & \rTo & \cC^0_\lm(A) & \rTo^\dif & \cC^1_\lm(A) & \rTo^\dif & \cC^2_\lm(A) & \rTo^\dif & \dots
\end{diagram} \]
\end{itemize}

For $n\geq 0$: the map $\imath:\sD^n\to\sE^n$ is inclusion;
the map $M:\sE^n\to\sF^n$ is defined to be
$(\sid-\bt)/2$ -- note that this differs from the map taken in \cite{Hel_CT} by a factor of~$2$;
and the map $N:\sF^n\to\sG^n$ is the averaging operator that was defined earlier.
For negative indices: we take $\imath:\sD^{-1}\to\sE^{-1}$ to be the identity map on $\cZ(A^*)$, and take $M:\sE^{-1}\to\sF^{-1}$, $N:\sF^{-1}\to\sG^{-1}$ to both be zero; in degrees $n\leq -2$ all the maps between complexes are necessarily zero. Then, recalling the identity \eqref{eq:intertwining}, it is clear that $\imath$ and $N$ are both chain maps. The proof that $M$ is a chain map is a straightforward calculation, which we omit.


\medskip
In order to apply the arguments of Section~\ref{s:abstract}, we need to define suitable horizontal splitting maps $j$, $P$ and $h$ in each degree.
For $n=-1$, this is trivial (see Figure~\ref{fig:stradivarius}), and for $n\leq -2$ we can take all maps to be zero. For $n\geq 0$, we take
$j=\imath:\cC^n_\lm(A)\to \cC^n(A)$, set $P=N:\cC^n_\lm(A)\to\cC^n(A)$,
and define $h:\cC^n(A)\to \cC^n(A)$ by
\[h = -\frac{2}{n+1}\sum_{k=1}^n k\bt^k\,. \]
%
%
Note that the norm of $h: \cC^n_\lm(A)\to \cC^n(A)$ is bounded above by $2(n+1)^{-1} \sum_{k=1}^n k  = n$\/.

\begin{figure}
\caption{Definitions in degrees $0$ and $-1$}
\label{fig:stradivarius}
\[\begin{diagram}[tight,width=3em,height=2.5em]
\cC^1_\lm(A) & \rTo^\imath & \cC^1(A) & \rTo^M & \cC^1(A) & \rTo^N & \cC^1_\lm(A) \\
\uTo & &  \uTo^{\dif}\dDots_{\sig_0} & & \uTo^{\dif'}\dDots_{\sig'_0} & & \uTo \\
\cC^0_\lm(A) & \rEq & \cC^0(A) & \rTo^0 & \cC^0(A) & \rEq & \cC^0_\lm(A) \\
\uTo^{\inc} & &  \uTo^{\inc} & & \uTo\dDots & & \uTo \\
 \cZ(A^*) & \rEq & \cZ(A^*) & \rTo & 0 & \rTo & 0 
\end{diagram}\]
\end{figure}

\begin{lem}\label{l:splitting}
Let $n\geq 0$\/.
Then $\imath P + hM = \sid$ and $Mh+jN=\sid$, regarded as maps $\cC^n(A)\to\cC^n(A)$\/.
\end{lem}

\begin{proof}
It suffices to prove the first identity, since $hM=Mh$ and $jN=\imath P$. This can be checked by direct calculation,~{\it viz}\/. 
\[ \begin{aligned}
\imath P + hM
 & = \frac{1}{n+1}\sum_{j=0}^n \bt^j -\frac{1}{n+1}\sum_{k=1}^n k\bt^k(\sid-\bt) \\
 & = \frac{1}{n+1}\left[ \sid + \sum_{j=1}^n \bt^j -\frac{1}{n+1}\sum_{k=1}^n k\bt^k + \sum_{j=2}^n (j-1)\bt^j + n\bt^{n+1} \right] \\
 & = \sid+\frac{1}{n+1}\left[ \sum_{j=1}^n \bt^j -\frac{1}{n+1}\sum_{k=1}^n k\bt^k + \sum_{j=2}^n (j-1)\bt^j \right] 
 & = \sid\,,
\end{aligned} \]
as required.
\end{proof}

\end{subsection}

\begin{subsection}{Traces}
Given $\psi\in A^*$ and $n\geq 0$\/, let $\psi^{(n)}\in \cC^n(A)$ be the cochain defined by
\[ \psi^{(n)}(a_1,\ldots, a_n)(a_0) \defeq \tau(a_0\dotsb a_n) \]
The following result is easily verified by a direct calculation, and we omit the proof.

\begin{lem}
If $\tau$ is a continuous trace on $A$\/, then $\tau^{(2n)}$ is a cyclic cocycle.
\end{lem}

We shall need to know how cyclic cocycles of this form transform under the shift map~$\widetilde{S}$.
Observe that for any $\psi\in A^*$\/, we have
$\dif'(\psi^{(2n)}) = \psi^{(2n+1)}$.
(This follows from a direct calculation, observing that the formula for $\dif'\tau^{(2n)}$ consists of $2n+1$ terms which cancel pairwise, save for the last one).
A similar calculation shows that $\dif(\tau^{(2n+1)})=\tau^{(2n+2)}$ if $\tau$ is a continuous trace on~$A$\/.

\begin{prop}
Let $\tau\in \cZ(A^*)$. Then $\widetilde{S}(\tau^{(2n)}) = S^\natural(\tau^{(2n)}) = \tau^{(2n+2)}$ for all $n\geq 0$\/.
\end{prop}

\begin{proof}
Since $\tau$ is a trace, $\bt (\tau^{(2n+1)}) = - \tau^{(2n+1)}$\/.
Hence 
\[ h\dif'(\tau^{(2n)})
 = h(\tau^{(2n+1)})
 = - \frac{2}{2n+2}\sum_{k=1}^{2n+1}(-1)^k k \tau^{(2n+1)}
 = \tau^{(2n+1)} 
 \,. \]
Therefore $S^\natural(\tau^{(2n)}) = \dif (\tau^{(2n+1)}) = \tau^{(2n+2)}$, and since this is already cyclic it is unchanged after we apply the averaging projection~$P$\/. 
\end{proof}

\end{subsection}

\begin{subsection}{The cyclic cohomology of biflat algebras}
We are almost ready to prove Theorem~\ref{t:newNEEDED}. In order to apply the results of Section~\ref{s:abstract}, we now have to impose some additional conditions on the complexes $\sE^\blob$ and $\sF^\blob$\/.

\paragraph{Condition 1.}
There exists a contracting homotopy $\sig$ for the complex
\[ \begin{diagram}[tight,width=2.8em,height=2em]
\dots\to & 0 & \rTo& \cZ(A^*) & \rTo^{\inc} & \cC^0(A) & \rTo^\dif & \cC^1(A) & \rTo^\dif & \cC^2(A) & \rTo^\dif &\dots
\end{diagram} \]

\paragraph{Condition 2.}
There exists a contracting homotopy $\sig'$ for the complex
\[ \begin{diagram}[tight,width=2.8em,height=2em]
\dots \to & 0 & \rTo& 0 & \rTo & \cC^0(A) & \rTo^{\dif'} & \cC^1(A) & \rTo^{\dif'} & \cC^2(A) & \rTo^{\dif'} & \dots
\end{diagram} \]

Moreover, to simplify some of the ensuing estimates, it is convenient to require that the following holds.

\paragraph{Condition 3.}
\[ c_1 \defeq \sup_n \norm{\sig_n: \cC^{n+1}(A) \to \cC^n(A)} \]
and
\[ c_2\defeq \max\left( 1, \sup_n\norm{\sig'_n: \cC^{n+1}(A) \to \cC^n(A)}\right) \]
are both finite.

\begin{lem}\label{l:biflat-marche}
Suppose $A$ is biflat with constant $K$. Then Conditions 1, 2 and~3 are satisfied, with $c_1=c_2=K$\/.
\end{lem}

\begin{proof}[Sketch of proof]
This is mostly standard, known material but with some minor additional book-keeping. (The fact that Condition~1 is satisfied can be found as \cite[Propo\-sition~2.8.62]{Dal_BAAC}, but for the reader's convenience we shall sketch the proof below.)
 
Let $\pi:A\ptp A\to A$ be the product map.
By hypothesis, there exists a bounded linear $A$-bimodule map $\rho:(A\ptp A)^*\to A^*$ which is left inverse to $\pi^*$ and has norm $\leq K$\/.
%
For $n\geq 1$, we define $\sig_{n-1}=\sig'_{n-1}: \cC^n(A)\to\cC^{n-1}(A)$ by
\begin{equation}\label{eqn:splitting_dfn}
\sig'_{n-1}\psi(a_0,\ldots,a_{n-1}) = \sig\psi(a_0,\ldots, a_{n-1}) \defeq \rho\left[\psi(\blank,\blank, a_1,\dots,a_{n-1}) \right](a_0)
\end{equation}
Then, for $n\geq 1$, $\psi\in\cC^n(A)$, and $a_0,\dots, a_n\in A$, a direct calculation of 
\[ \dif\sig_{n-1}(\psi)(a_0,\dots,a_n) + \sig_n\dif(\psi)(a_0,\dots,a_n)\]
shows that most terms cancel to leave us with
\[ \begin{aligned}
	  & \rho\left[ \psi(\blank,\blank,a_2,\dots,a_n) \right](a_0a_1) \\
	+ (-1)^n & \rho\left[ \psi(\blank,\blank,a_1,\dots,a_{n-1}) \right](a_na_0) \\
	+ & \rho\left[ \psi(\blank\cdot\blank,a_1,\dots, a_n) \right](a_0) \\
	- & \rho\left[ \psi(\blank,\blank\cdot a_1,a_2,\dots,a_n) \right](a_0) \\
	+ (-1)^{n-1} & \rho\left[ \psi(a_n\cdot\blank,\blank, a_1,\dots, a_{n-1}) \right](a_0)\,;
\end{aligned} \]
and, since $\rho$ is an $A$-bimodule map and is left inverse to~$\pi$, this in turn reduces to $\psi(a_0,\dots,a_n)$. Hence
\[ \dif\sig_{n-1}(\psi)+\sig_n\dif(\psi) = \psi \qquad\text{ for all $\psi\in\cC^n(A)$, $n\geq 1$\/.} \]
A similar calculation shows that $\dif'\sig'_{n-1}(\psi)+\sig'_n\dif'\psi=\psi$ for all such $\psi$. 

It only remains to define $\sig$ and $\sig'$ appropriately in negative degrees. Since $\sig'_0\dif'=\sid$, as shown by a quick calculation, we can put $\sig'_n=0$ for all $n\leq -1$ and satisfy Condition~2.
Now observe that, since $\sig_1\dif+\dif\sig_0=\sid$\/, we have
\[ \dif\sig_0\dif = (\dif\sig_0+\sig_1\dif)\dif = \dif\/; \]
therefore, on putting $\sig_{-1}=\sid-\sig_0\dif:\cC^0(A)\to \cZ(A^*)$, we find that $\dif\sig_{-1}+\sig_0\dif$ is the identity map on $\cC^0(A)$ and $\sig_{-1}\inc$ is the identity map on $\cZ(A^*)$. Thus, putting $\sig_n=0$ for all $n\leq -2$, we have satisfied Condition~1. Finally, it is clear from our construction that Condition~3 is satisfied with $c_1=c_2=\norm{\rho}$.
\end{proof}

\begin{rem}\label{rem:cyclic-der}
By our assumption on the splitting homotopy $\sig$\/, if $\psi\in \cC^1(A)$ then $\psi=\sig\dif(\psi)+\dif\sig(\psi)$\/. In particular, if $\psi$ is a derivation then $\psi=\dif\sig(\psi)$ is inner, and thus cyclic.
\end{rem}

We suppose for the rest of this section that $A$ is a Banach algebra satisfying Conditions 1, 2 and~3. The first two of these conditions are known to imply (by assembling the appropriate results from~\cite{Hel_CT}) that the shift map $S:\cH_\lm^n(A)\to\cH_\lm^{n+2}$ is an isomorphism for each~$n$; the extra constraints imposed by Condition~3 allow us to give the following quantitative version.

\begin{prop}\label{p:ind-step}
Let $n\geq 0$. For every $\psi\in \cZ^{n+2}_\lm(A)$\/, there exists $\chi\in \cC^{n+1}_\lm(A)$ and $\varphi\in \cZ^n_\lm(A)$\/, such that
\[ \psi= \widetilde{S}\varphi + \dif\chi\quad\text{ and }\quad \norm{\chi} \leq (n+1)^2c_1 c_2 \norm{\psi} \quad\text{ and }\quad \norm{\varphi}\leq c_1 c_2 \norm{\psi}\,.\]
\end{prop} 

\begin{proof}
Put $\varphi \defeq \widetilde{R}(\psi)$ and $\chi= -T(\psi)$\/. Then since
\[ \widetilde{S}\widetilde{R}\psi - \psi = (\dif T+T\dif)(\psi) = \dif T\psi\,,\]
we have $\psi=\widetilde{S}\varphi+\dif\chi$\/.
The norm estimates follow, since $\norm{\widetilde{R}}\leq \norm{N\sig'M\sig\imath} \leq c_1 c_2$, while
$\norm{T} \leq \norm{T^\natural} \leq \norm{\dif h \sig' M\sig}+\norm{\sig}
   \leq (n+2)\cdot n\norm{\sig'}\cdot \norm{\sig} + \norm{\sig}
   \leq (n+1)^2 c_2 c_1$\/.
\end{proof}

In view of Lemma~\ref{l:biflat-marche}, Theorem~\ref{t:newNEEDED} will now follow from the following result (which is slightly more precise).

\begin{thm}\label{t:NEEDED}
Let $m\geq 0$\/.
\begin{YCnum}
\item For every $\psi\in \cZ^{2m+1}_\lm(A)$, there exists $\chi\in \cC^{2m}_\lm(A)$ such that $\psi=\dif\chi$ and
\[ \norm{\chi}\leq 2(m+1)^3(c_1c_2)^{2m}\norm{\psi} \,.\]   
\item For every $\psi\in \cZ^{2m+2}_\lm(A)$, there exist $\chi\in \cC^{2m+1}_\lm(A)$
 and $\tau\in \cC^0_\lm(A)$
 such that
$\psi=\tau^{(2m+2)}+\dif\chi$\/, and which satisfy
\[ \norm{\tau}\leq (c_1c_2)^{m+1}\norm{\psi} \quad,\quad \norm{\chi}\leq 2(m+1)^3(c_1c_2)^{2m+1}\norm{\psi}\,.\]
\end{YCnum}
\end{thm}

\begin{proof}[Proof of Theorem~\ref{t:NEEDED}]
We do each part by induction. For cochains in even degree: the case $m=0$ is given by the preceding proposition with $n=0$\/.
Suppose it holds for $m=k-1$, where $k\geq1$\/, and let $\psi\in \cZ^{2k+2}_\lm(A)$. Applying the proposition with $n=2k$\/, there exists $\chi\in \cC^{2k+1}_\lm(A)$
and $\varphi\in \cZ^{2k}_\lm(A)$
such that
\[ \psi= \dif\chi + \widetilde{S}\varphi \quad,\quad\norm{\chi}\leq (2k+1)^2c_1c_2\norm{\psi}\quad,\quad\text{ and } 
\norm{\varphi}\leq c_1c_2\norm{\psi}\,.\]
 
Applying the inductive hypothesis to $\varphi$ yields $\tau\in \cC^0_\lm(A)$ and $\chi_2\in \cC^{2k-1}_\lm(A)$ such that
\[ \varphi=\tau^{(2k)}+\dif\chi_2 \quad\text{ and }\quad
  \norm{\tau}\leq (c_1c_2)^k \norm{\varphi} \quad\text{ and }\quad
  \norm{\chi_2}\leq 2k^3(c_1c_2)^{2k-1} \norm{\varphi}\,. \]
Recalling that $\widetilde{S}\dif+\dif\widetilde{S}=0$\/, we therefore have
\[
\psi
  = \dif\chi - \dif\widetilde{S}(\chi_2) + \widetilde{S}\tau^{(2k)} 
  = \dif(\chi-\widetilde{S}\chi_2) +\tau^{(2k+2)}\;,
\]
where
\[ \norm{\tau} \leq (c_1c_2)^k\norm{\varphi} \leq (c_1c_2)^{k+1}\norm{\psi}\,,\]
and
\[ \begin{aligned}
\norm{\chi-\widetilde{S}\chi_2}
 & \leq (2k+1)^2c_1c_2\norm{\psi} + \norm{\widetilde{S}} \cdot 2k^3 (c_1c_2)^{2k-1} \norm{\varphi} \\
 & \leq (2k+1)^2(c_1c_2)^{2k+1} \norm{\psi} + c_1c_2 \cdot 2k^3 (c_1c_2)^{2k-1} \cdot c_1c_2\norm{\psi} \\
 & \leq 2(k+1)^3 (c_1c_2)^{2k+1} \norm{\psi}\,,
\end{aligned} \]
and this completes the inductive step.

The proof for cochains in odd degree is similar. 
For $m=0$\/, the claim follows from Remark~\ref{rem:cyclic-der}.
If the claim holds for $m=k-1$ where $k\geq1$\/, let $\psi\in \cZ^{2k+1}_\lm(A)$\/. Applying Proposition~\ref{p:ind-step} with $n=2k-1$\/, we obtain $\chi\in\cC^{2k}_\lm(A)$ and $\varphi\in\cZ^{2k-1}_\lm(A)$ such that
\[ \psi= \dif\chi + \widetilde{S}\varphi \quad,\quad\norm{\chi}\leq (2k)^2c_1c_2\norm{\psi}\quad,\quad\text{ and } 
\norm{\varphi}\leq c_1c_2\norm{\psi}\,.\]
Applying the inductive hypothesis to $\varphi$ yields $\chi_2\in \cC^{2k-2}_\lm(A)$ such that
\[ \varphi=\dif\chi_2 \quad\text{ and }\quad \norm{\chi_2} \leq 2k^3(c_1c_2)^{2k-2}\norm{\varphi}\,. \] 
Recalling that $\widetilde{S}\dif+\dif\widetilde{S}=0$\/, we therefore have
\[ \psi= \dif\chi +\widetilde{S}\dif\chi_2 = \dif(\chi -\widetilde{S}\chi_2) \,,\]
where
\[ \begin{aligned}
\norm{\chi-\widetilde{S}\chi_2}
  \leq \norm{\chi} + \norm{\widetilde{S}}\norm{\chi_2} 
 & \leq (2k)^2c_1c_2 \norm{\psi} + c_1c_2\cdot 2k^3 (c_1c_2)^{2k-2}\cdot c_1c_2\norm{\psi} \\
 & \leq 2(k+1)^3 (c_1c_2)^{2k} \norm{\psi}\,,
\end{aligned} \]
completing the inductive step.
\end{proof}

\end{subsection}

\end{section}

\section*{Acknowledgments}
Some of the work in the preprint \cite{CGW_bands}, which motivated the present article, was presented at the 19$^\mathrm{th}$ International Conference on Banach Algebras (B\c{e}dlewo, July 14--24,~2009).
The author thanks the Polish Academy of Sciences, the
European Science Foundation under the ESF-EMS-ERCOM partnership, and the
Faculty of Mathematics and Computer Science of the Adam Mickiewicz University
at Pozna\'n for their support of this meeting.
He would also like to thank the referee of this article, for bringing \cite{Quillen_IHES88} to his attention, and for pointing out that it discusses the cyclic cochain complex and its cohomology in some depth.

\appendix


\begin{section}{The SBI sequence, revisited}\label{s:appendix}
In this appendix, we will construct a long exact sequence relating the cohomology of the complexes $\sD^\blob$, $\sE^\blob$, $\sF^\blob$ and $\sG^\blob$, under certain additional hypotheses on the complex $\sF^\blob$. When specialized to the setting of Section~\ref{s:appl}, we will recover a special case of the Connes-Tsygan exact sequence for Banach algebras as constructed in~\cite{Hel_CT}. The results thus obtained are less general than those of \cite{Hel_CT}, but still apply to Banach algebras that are one-sided flat over themselves, for example.

Our approach here is direct: this avoids the need to discuss various general principles and connecting homomorphisms from homological algebra, and keeps our discussion more self-contained, but means that the motivation\footnotemark\
for some of our calculations is rather obscure. A slightly more conceptual approach can be found in the standard sources \cite[Part II]{Connes_IHES_NCDG} and 
\cite{Hel_CT}, albeit in the particular setting of the cyclic cochain complex.
\footnotetext{%
Readers well versed in homological algebra will recognize that we are tracing through the construction, by hand, of a so-called ``long exact sequence of Gysin type'' -- such a sequence is often derived from more general arguments with spectral sequences, but it seems worth recording how to do this directly, for sake of convenient future reference.
}

\paragraph{The insertion map.}
The inclusion $\imath: \sD^\blob\to \sE^\blob$ induces a morphism on cohomology, which we denote by $I:H^n(\sD^\blob) \to H^n(\sE^\blob)$. In particular, we have a diagram
\[ \begin{diagram}[tight,width=4em]
H^n(\sG^\blob) & \rTo^S & H^{n+2}(\sD^\blob) & \rTo^I & H^{n+2}(\sE^\blob)
\end{diagram} \]
for each~$n$\/.
\begin{lem}\label{l:starter}
$\ker I \subseteq \image S$.
\end{lem}
\begin{proof}
Since
\[ \begin{aligned}
\imath \widetilde{S}
 &  = \imath P \delta h\delta' j 
 & = (\sid - hM) \delta h \delta' j
\end{aligned} \]
we have
\[ \begin{aligned}
\imath \widetilde{S} -  \delta h \delta'j
 & = - hM\delta h \delta' j \\
 & = h \dif' j \dif'' N j & \quad\text{(by Lemma~\ref{l:useful1}(ii))} \\
 & = \;\; h\delta'j\delta'' & \quad\text{($N$ is a chain map and $Nj=\sid$).} \\
\end{aligned} \]
Thus $\imath\widetilde{S}(Z^n(\sG^\blob))\subseteq \dif'(\sE^{n+1})$, and so on passing to cohomology we have $IS=0$\/.
\end{proof}

\medskip
\paragraph{Constructing $\widetilde{B}$.}
From here onwards, we suppose that there exists a splitting homotopy $\sig': \sF^\blob\to\sF^{\blob-1}$ for the complex~$\sF^\blob$\/.

Define $\widetilde{B}: \sE^{n+1} \to \sG^n$ to be the composite mapping.
\[ \begin{diagram}
 \sE^{n+1} & \rTo^M & \sF^{n+1} & \rTo^{\sig'} & \sF^n & \rTo^N & \sG^n 
\end{diagram} \]

\begin{prop}[{see \cite[Part II, Lemma 30]{Connes_IHES_NCDG}}]
\label{p:anticom}
$\dif'' \widetilde{B}+\widetilde{B}\dif =0$\/.
\end{prop}

\begin{proof}
This is immediate from Lemma~\ref{l:useful2}.
\end{proof}

It follows that $\widetilde{B}$ descends to a well-defined map of cohomology, which we denote by $B: H^{n+1}(\sE^\blob) \to H^n(\sG^\blob)$ for each~$n$\/. This fits together with the maps $I$ and $S$ that we have already defined, to yield a sequence of maps

\begin{equation}\label{eq:SBI}
\begin{diagram}[tight,width=3.5em]
\dots  & H^{n-1}(\sG^\blob) & \rTo^S & H^{n+1}(\sD^\blob) & \rTo^I & H^{n+1}(\sE^\blob) & \rTo^B & H^n(\sG^\blob) & \rTo^S &\dots \tag{$\dagger$}
\end{diagram}
\end{equation}

\begin{thm}[Abstract version of the Connes-Tsygan exact sequence]
\label{t:SBI}
The sequence \eqref{eq:SBI} is exact.
\end{thm}

%


The proof will be broken up into several steps.

\begin{prop}\label{p:BI}
$\ker B=\image I$\/.
\end{prop}

\begin{proof}
Since $\widetilde{B}\imath=0$\/, it is immediate that $\ker B \supseteq \image I$\/. Conversely, let $\psi\in Z^{n+1}(\sE^\blob)$ be such that $\widetilde{B}(\psi)=\dif''\varphi$ for some $\varphi\in \sG^{n-1}$. We seek $\chi\in\sE^n$ such that $\psi-\dif\chi \in \image \imath$\/; equivalently, since $\image \imath=\ker M$, it is enough to find $\chi\in\sE^n$ such that $M(\psi)=M\dif(\chi)=\dif'M(\chi)$.

We start by noting that, since $\dif\psi=0$\/,
\begin{equation}\label{eq:larry}
\dif'\sig' M(\psi) = (\sid-\sig'\dif')M(\psi) = M(\psi) - \sig' M\dif(\psi) = M(\psi).
\end{equation}
Now
\begin{equation}\label{eq:curly}
\begin{aligned}
\dif'\sig'M(\psi) = \dif'(jN+Mh)\sig'M(\psi)
 & = \dif'j\widetilde{B}(\psi) + \dif'Mh\sig'M(\psi) \\
 & = \dif'j\dif''(\varphi) + \dif'M (\chi_2)
\end{aligned}
\end{equation}
where we have put $\chi_2\defeq h\sig'M(\psi)$.
Next, since $N\dif'j =Nj\dif''=\dif''$, we see that
\begin{equation}\label{eq:moe}
\dif'j\dif'' = \dif'jN\dif' j = \dif'(\sid-Mh)\dif'j = -\dif'Mh\dif' j\,.
\end{equation}
Thus, if we put $\chi_1 = -h\dif'j (\psi)$, combining \eqref{eq:curly} and \eqref{eq:moe} gives
\[ \dif'\sig'M\psi = \dif'M(\chi_1) + \dif'M(\chi_2) \]
and combining this with \eqref{eq:larry} yields $M(\psi) = \dif'M(\chi_1+\chi_2)$ as required.
\end{proof}

We next prove that $\ker S=\image B$. This requires some preparatory lemmas.

\begin{lem}\label{l:ghastly}
\[
 \widetilde{S}\widetilde{B} = P\dif - \dif P  - P\dif h \sig' M\dif + \dif P\dif h\sig' M
\]
\end{lem}

\begin{proof}[Proof of Lemma~\ref{l:ghastly}]
We have
\begin{flalign*}
\quad \widetilde{S}N  = P\dif h\dif'j N 
  = P\dif h\dif'(\sid-Mh) 
 & = P\dif h\dif' - P\dif h\dif' Mh \\
 & = P\dif h\dif' + P\imath\dif P\dif h & \text{(by Lemma~\ref{l:useful1}(i))}\quad \\
 & = P\dif h\dif' + \dif P\dif h & \text{(since $P\imath=\sid$).}\quad
\end{flalign*}
Thus $\widetilde{S}N-\dif P\dif h = P\dif h\dif'$. Since $\widetilde{B}=N\sig'M$, this implies that
\[
\begin{aligned}
 \widetilde{S}\widetilde{B} - \dif P\dif h\sig' M
  & = P\dif h\dif'\sig'M \\
  & = P\dif h(\sid-\sig'\dif') M \\
  & = P\dif hM - P\dif h \sig' M\dif & \quad\text{(since $M$ is a chain map)} \\
  & = P\dif (\sid - \imath P) - P\dif h \sig' M\dif & \\
  & = P\dif - \dif P  - P\dif h \sig' M\dif & \quad\text{(since $\imath$ is a chain map,} \\
 & & \quad\text{and $P\imath=\sid$\/.)} 
\end{aligned}
\]
This completes the proof.
\end{proof}

Define $Y: \sG^n \to E^{n+1}$ to be the composite map
\[ \begin{diagram}
\sG^n & \rTo^j & \sF^n & \rTo^{\dif'} & \sF^{n+1} &\rTo^h & \sE^{n+1}
\end{diagram} \]
In general there is no reason for $Y$ to be a chain map. However, we do have the following useful identity.

\begin{lem}\label{l:Dazzler}
$\widetilde{B}Y = \sid - \dif''N\sig' j - N\sig'j\dif''$.
\end{lem}
\begin{proof}
We have
\begin{equation}
\label{eq:Frankie}
\begin{aligned}
\widetilde{B}Y = N\sig'M\cdot h\dif' j
 & = N\sig'(\sid-jN)\dif' j \\
 & = N\sig'\dif' j - N\sig'jN\dif'j \\
 & = N\sig'\dif' j - N\sig'j\dif''Nj & \text{(since $N$ is a chain map)}\quad \\
 & = N\sig'\dif' j - N\sig'j\dif'' & \text{(since $Nj=\sid$),}\quad \\
\end{aligned}
\end{equation}
while
\begin{equation}\label{eq:Hugh}
\begin{aligned}
N\sig'\dif' j = N(\sid-\dif'\sig')j 
 & = \sid - N\dif'\sig'j & \quad\text{(since $Nj=\sid$)} \\
 & = \sid - \dif''N\sig' j & \quad\text{(since $N$ is a chain map).}
\end{aligned}
\end{equation}
Combining \eqref{eq:Frankie} and \eqref{eq:Hugh} concludes the proof.
\end{proof}

Note that $\widetilde{S}=P\dif Y$.

\begin{prop}\label{p:SB}
$\ker S = \image B$.
\end{prop}

\begin{proof}
It follows from Lemma~\ref{l:ghastly} that $SB=0$, i.e.~that $\ker S \supseteq \image B$. The converse inclusion is proved as follows.
Let $\psi\in Z^n(\sG^\blob)$ be such that $S([\psi])=0$. Then there exists $\varphi\in \sD^{n+1}$ such that $\widetilde{S}(\psi)=\dif\varphi$\/. Consider $Y(\psi)$: although this might not lie in $Z^{n+1}(\sE^\blob)$, we have $P\dif Y(\psi) = \dif\varphi$\/, and so 
\[ \begin{aligned}
\dif\imath(\varphi) =\imath\dif(\varphi)
   = \imath P \dif(\chi) 
   = (\sid - hM)\dif(\chi) 
   = \dif(\chi) -hM\dif(\chi).
\end{aligned} \]
Thus
\begin{flalign*}
\qquad \dif(\imath\varphi - Y\psi)
  = - hM\dif Y(\psi)  
 & = - hM\dif h\dif' j(\psi) \\
 & = \;\; h\dif'j\dif''N j(\psi) & \text{(by Lemma~\ref{l:useful1}(ii))}\quad  \\
 & = \;\; h\dif'j\dif''(\psi) & \text{(since $Nj=\sid$)} \quad \\
 & = \;\; 0  & \text{(since $\psi\in Z^n(\sG^\blob)$).}\quad 
\end{flalign*}
Putting $\chi\defeq \imath\varphi -Y(\psi)$, we therefore have $\chi \in Z^{n+1}(\sE^\blob)$. Now since $\widetilde{B}\dif= -\dif\widetilde{B}$, we see that $\widetilde{B}(\chi)\in Z^n(\sG^\blob)$; while
\[ \begin{aligned}
\widetilde{B}(\chi)
 & = -\widetilde{B}Y(\psi) & \quad\text{(since $\widetilde{B}\imath=0$)} \\ 
 & = -\psi + \dif''N\sig' j(\psi)  & \quad\text{(by Lemma~\ref{l:Dazzler}).}
\end{aligned} \]
Hence $B([-\chi])= [\psi]$, and thus $\ker S \subseteq \image B$\/.
\end{proof}

\begin{prop}\label{p:IS}
$\ker I = \image S$.
\end{prop}

\begin{proof}
By Lemma~\ref{l:starter}, we have $\ker I\supseteq \image S$. Conversely, let $\psi\in Z^{n+2}(\sD^\blob)$ and suppose that $\imath(\psi)\in B^{n+2}(\sE^{\blob})$\/. Then $\imath(\psi)=\delta(\varphi)$ for some $\varphi\in \sE^{n+1}$\/.
Put $\chi \defeq \widetilde{B}(\varphi)\in \sG^n$\/, and note that
\[ \begin{aligned}
\dif''\chi = \dif''\widetilde{B}\varphi 
 & = -\widetilde{B}\dif\varphi	& = -\widetilde{B}\imath(\psi)  & = 0 \qquad\\
 & \text{(by Propn \ref{p:anticom})} & & \text{(since $\widetilde{B}\imath=0$).}
\end{aligned} \]
Thus $\chi\in Z^n(\sG^\blob)$. We claim that $\widetilde{S}(\chi)$ is cohomologous to $\psi$\/, which will show that $S([\chi])=[\psi]$, and hence that $\ker I=\image S$\/.
The claim follows by applying Lemma~\ref{l:ghastly}, {\it viz}\/.
\[ \begin{aligned}
\widetilde{S}(\chi) =
\widetilde{S}\widetilde{B}(\varphi)
 & = (\dif P\dif h\sig' M + P\dif - \dif P - P\dif h\sig' M\dif)(\varphi) \\
 & = \dif P\dif h\sig' M(\varphi) + \psi- \dif P(\varphi),
\end{aligned} \]
where we use the identities $P\dif(\varphi)=P\imath(\psi)=\psi$\/, and $M\dif(\varphi) = M\imath(\psi)=0$\/.
\end{proof}

Combining Propositions \ref{p:BI}, \ref{p:SB} and \ref{p:IS}, we conclude that the diagram~\eqref{eq:SBI} is an exact sequence, and Theorem~\ref{t:SBI} is proved.

\end{section}


\begin{thebibliography}{1}

\bibitem{CGW_bands}
Y. Choi, F. Gourdeau, M. C. White, {\em Simplicial cohomology of band semigroup algebras}, preprint. See arXiv \href{http://arxiv.org/abs/1004.2301}{\tt math.1004.2301}.

\bibitem{Connes_IHES_NCDG}
A. Connes, {\em Noncommutative differential geometry}, Inst. Hautes \'Etudes Sci. Publ. Math., 62 (1985), 41--144.

\bibitem{Dal_BAAC}
H. G. Dales, {\em {B}anach algebras and automatic continuity}, vol.~24 of London Mathematical Society Monographs (New Series), Clarendon Press, Oxford University Press, New York, 2000.

\bibitem{Hel_HBTA}
A.~{\relax Ya}. Helemski{\u\i},
{\em The {H}omology of {B}anach and {T}opological {A}lgebras}, vol.~41 of Mathematics and its Applications (Soviet Series), Kluwer Academic Publishers Group, Dordrecht, 1989.

\bibitem{Hel_CT}
A.~{\relax Ya}. Helemski{\u\i},
{\em Banach cyclic cohomology and the {C}onnes-{T}sygan exact sequence}, J.~London Math Soc. (2), 46 (1992),  449--462.

\bibitem{Lod_CH_ed1}
J.-L. Loday, {\em Cyclic homology}, vol.~301 of Grundlehren der
  Mathematischen Wissen\-schaften,
 Springer-Verlag, Berlin, 1992.
\newblock Appendix E by Mar\'\i a O. Ronco.

\bibitem{Quillen_IHES88}
D. Quillen, {\em Algebra cochains and cyclic cohomology}, Inst. Hautes \'Etudes Sci. Publ. Math., 68 (1988), 139--174.

\end{thebibliography}

\vfill\eject

\vfill

\noindent%
\begin{tabular}{l}
Y. Choi\\
D\'epartement de math\'ematiques
et de statistique,\\
Pavillon Alexandre-Vachon\\
Universit\'e Laval\\
Qu\'ebec, QC \\
Canada, G1V 0A6 \\
{\bf Email: \tt y.choi.97@cantab.net}
\end{tabular}

\end{document}